\documentclass[11pt]{amsart}

\usepackage{amsmath,amssymb}
\usepackage{fullpage}
\usepackage{amscd}
\usepackage{amsfonts, amssymb, stmaryrd}
\usepackage[all]{xy}
\usepackage{tikz}
\usepackage{latexsym}
\usepackage{cite}
\usepackage{url}
\usepackage{boxedminipage}
\usepackage{amsbsy}

\usepackage{hyperref}
\hypersetup{colorlinks}

\usepackage{enumitem}

\usepackage{dutchcal}

\newcommand{\fg}{\mathfrak g}

\newcommand{\cA}{{\mathcal A}}

\newcommand{\cah}{\mathcal{h}}
\newcommand{\cp}{\mathcal{p}}

\newcommand{\cM}{{\mathcal M}}

\newcommand{\N}{\mathbb{N}}
\newcommand{\R}{\mathbb{R}}

\newcommand{\Z}{\mathbb{Z}}

\usepackage{mathrsfs}

\newcommand{\scA}{\mathscr{Ac}}
\newcommand{\scB}{\mathscr{B}}
\newcommand{\scC}{\mathscr{C}}
\newcommand{\scI}{\mathscr{I}}

\newcommand{\scV}{\mathscr{V}}

\newcommand{\DiffSp}{\mathsf{DiffSp}}

\DeclareMathAlphabet{\mathpzc}{OT1}{pzc}{m}{it}

\newcommand{\inv}{^{-1}}

\DeclareMathOperator{\supp}{supp}

\newcommand{\cin}{C^\infty}

\numberwithin{equation}{section}

\theoremstyle{definition}

\newtheorem{lemma}[equation]{Lemma}
\newtheorem{theorem}[equation]{Theorem}

\newtheorem{corollary}[equation]{Corollary}
\newtheorem*{corollary*}{Corollary}

\newtheorem{definition}[equation]{Definition}
\newtheorem{remark}[equation]{Remark}
\newtheorem{example}[equation]{Example}
\newtheorem{notation}[equation]{Notation}

\begin{document}
\title{Stability and bifurcations of symmetric tops}
\author{Eugene Lerman}

\setcounter{tocdepth}{1}

\begin{abstract}  
We study the stability and bifurcation of relative equilibria of a particle on the
Lie group $SO(3)$ whose motion is governed by an $SO(3)\times SO(2)$
invariant metric and an $SO(2)\times SO(2)$ invariant potential.  Our
method is to reduce the number of degrees of freedom at {\em singular}
values of the $SO(2)\times SO(2)$ momentum map and study the stability
of the equilibria of the reduced systems as a function of spin.  The
result is an elementary analysis of the fast/slow transition in the
Lagrange and Kirchhoff tops.

More generally, since an $SO(2)\times
SO(2)$ invariant potential on $SO(3)$ can be thought of as ${\Z}_2$
invariant function on a circle,  we analyze the stability and
bifurcation of relative equilibria of the system in terms of the
second and fourth derivative of the function.
  \end{abstract}
\maketitle
\tableofcontents

\section{Introduction}

Recall the geometric mechanics approach to classical systems which was
developed in the early 1960s: A ``simple'' classical mechanical system consists of a
manifold $Q$, the configuration space of the system, a Riemannian metric $g$ in $Q$,
the kinetic energy, 
and a function $V:Q\to \R$, the potential.   These data define a
Hamiltonian $H:T^*Q\to \R$ on the cotangent bundle of $Q$ which 
 is given by
\[
H(q, p) = \frac{1}{2} g^*_q(p,p) +V(q)
\]
where $q\in Q$, $p\in T^*_qQ$ and $g^*\in \mathrm{Sym}^2(T^*Q)$ is the
dual metric.  The Hamiltonian $H$ together with the canonical
symplectic form $\omega_{T^*Q}$ on the cotangent bundle give rise to a
vector field $\Xi_H:T^*Q\to T(T^*Q)$ which is uniquely defined by the equation
\[
  \omega_{T^*Q}(\Xi_H, \cdot) = dH.
\]
An action of a Lie group $G$ on $Q$ that preserves the metric $g$ and
the potential $V$ lifts to an action of $G$ on the cotangent bundle
$T^*Q$ that preserves the Hamiltonian $H$ and the symplectic form
$\omega_{T^*Q}$.   Noether's theorem in this setting translates into the existence of an
equivariant moment map $\mu:T^*Q \to \fg^*$ ($\fg^*$ denotes the dual
of the Lie algebra $\fg$ of the Lie group $G$) with the property that
$\mu$ is constant along the intergral curves of the vector field
$\Xi_H$.   Since the Hamiltonian $H$ and the symplectic form are
$G$-invariant the Hamiltonian vector field $\Xi_H$ is $G$-invariant as
well.  Hence its flow $\Phi^H_t$ is $G$-equivariant.  Noether's
theorem combined with equivariance of the flow implies that we have
continuous time dynamical systems on the topological spaces
$\{\mu\inv(\alpha)/G_\alpha\}_{\alpha\in \fg^*}$, where $G_\alpha$
denotes the stabilizer of $\alpha\in \fg^*$ under the coadjoint
action.   When $\alpha$ is a regular value  of $\mu$ and the action
of $G$ is proper then, thanks to a theorem of Meyer \cite{Meyer} and of  Marsden and
Weinstein \cite{MW}, the topological space $M//_\alpha G:=
\mu\inv(\alpha)/G_\alpha$ is naturally a symplectic orbifold and the
flow induced by $\Phi^H_t$ is a flow of a Hamiltonian vector field.
The  theorem of Meyer and Marsden-Weinstein is known as  regular
symplectic  reduction and as (Meyer-)Marsden-Weinstein reduction (a
number of people seem to be unaware of Meyer's paper).\\

A {\sf top} is a classical mechanical system of the form
$(T^*SO(3), H(q,p) = \frac{1}{2} g^*_q(p,p) +V(q))$ where $SO(3)$ is
the special orthogonal group and $g$ is a left-invariant metric on the
group.  A top is {\sf symmetric} if two of its principal moments of
inertia are equal and the potential is invariant under the additional
symmetry.  This amounts to the metric $g$ being invariant
under the multiplicaiton  on the right by $SO(2)$ (here $SO(2)$ is the
subgroup of $SO(3)$ fixing the third standard basis vector
$e_3=(0,0,1)$), and the potential $V\in \cin(SO(3))$ being
$SO(2)\times SO(2)$-invariant. Here and elsewhere $SO(2)\times SO(2)$
acts on $SO(3)$ by multiplication on the left and right,
respectively. 

In a symmetric top one can view the rate of spin about its axis of
symmetry as a bifurcation parameter. Geometrically this amounts to
viewing a top as a family of Hamiltonian systems on the unit 2-sphere
$S^2$:
\[
\left\{ (T^*S^2, \omega_{T^*S^2} +r \omega_{S^2}, h(q,p) = \frac{1}{2}
g^*_q(p,p) +V(q)) \right\}_{r\in \R},
\]  
where $ \omega_{T^*S^2}$ is the canonical symplectic form on the
cotangent bundle,  $\omega_{S^2}$ is the standard area form on the
sphere, $g$ is the  standard round metric, and we have identified
the potential of the top with an $SO(2)$-invariant function on $S^2$.
See Corollary~\ref{cor:2.2} below.  Note that there is no loss of
generality in assuming that $r\geq 0$ since $r<0$ corresponds to spin
in the opposite direction.

Theorem~2 of \cite{Sch} implies that the algebra $\cin(S^2)^{SO(2)}$
of invariant functions on the sphere is isomorphic to
$\cin(S^1)^{\Z/2}$, where $\Z/2 =\{\pm 1\}$ acts on $S^1= \{(x,z)\in
\R^2 \mid x^2 + z^2 =1\}$ by $(-1)\cdot (x,z) = (-x, z)$.
Parameterize the upper half of $S^1$ by
\[
f:(-1, 1) \to S^1, \qquad f(u) = (u, \sqrt{1-u^2}).
\]
Then $f$ pulls back $\Z/2$-invariant functions on $S^1$ to $\Z/2$
invariant functions on $(-1,1)$ (where $-1\in \Z/2$ acts on $(-1,1)$
by $(-1) \cdot u = -u$).  It follows from a theorem of Whitney that
$\cin((-1, 1))^{\Z/2}$ is isomorphic to $\cin([0,1))$. That is, for
any $k\in \cin((-1, 1))^{\Z/2}$ there is a unique $\ell \in
\cin([0,1))$ with $k(u) = \ell (u^2)$.  (Recall that $\ell\in
\cin([0,1))$ iff there is $\varepsilon >0 $ and $\tilde{\ell}\in
\cin((-\varepsilon, 1))$ with $\ell = \tilde{\ell}|_{[0,1)}$.)
Therefore given a function $V\in \cin(S^2)^{SO(2)}$ there is a unique
function $W\in \cin([0,1))$ so that
\[
W(u^2) = V(u, 0, \sqrt{1-u^2}) \qquad \textrm{ for all } u\in [0,1). 
 \] 
We are now in position to formulate our stability and bifurcation
result.\\

\begin{theorem} \label{thm:1.1} Consider a 1-parameter family of $SO(2)$-invariant
Hamiltonian systems
\begin{equation}\label{eq:*}
\left\{ (T^*S^2, \omega_{T^*S^2} +r \omega_{S^2}, h(q,p) = \frac{1}{2}
g^*_q(p,p) +V(q)) \right\}_{r\geq 0},
\end{equation}
where $SO(2)$ acts by the lift of rotations about $e_3 =(0,0,1)$, $ \omega_{T^*S^2}$ is the canonical symplectic form on the
cotangent bundle,  $\omega_{S^2}$ is the standard area form on the
sphere, $g$ is the standard round round metric and $V \in
\cin(S^2)^{SO(2)}$.  Let $W\in \cin([0,1))$ be the function with
\[
W(u^2) = V(u, 0, \sqrt{1-u^2})
\]
for all $u\in [0,1)$.  Then
\begin{itemize}[leftmargin=0.2in]
\item[(i)] If $W'(0) >0$  then the straight up top (i.e. the point $(e_3,
  0)\in T^*S^2$ where $e_3 =(0,0,1)$) is stable for all values of $r$.
\item[(ii)]  Suppose  $W'(0)<0$ and $W''(0) > W'(0)$.  Then for $r> r_0 :=
\sqrt{-8W'(0)}$ the straight up top is stable.  As $r$ decreases
  below $r_0$ the top loses stability and we get a branch of stable
  relative equilibria bifurcating off the straight up position. That
  is, a ``Hamiltonian Hopf'' bifurcation takes place.
  \begin{tikzpicture}
    \draw[thick] (0,0) plot[domain=0:1.4] (\x,{2-(\x)^2}) node[anchor=south west] {};
 \draw[->] (0,0) -- (3,0) node[right] {$u$};
 \draw[thick, dashed] (0,0) -- (0,2) node[left]{$r_0$};   
 \draw[thick,->] (0,2) -- (0,3) node[above] {$r$};
\end{tikzpicture}

  \item[(iii)] Suppose $W'(0)<0$ and $W''(0)< W'(0)$.   Then for $0\leq r
  \leq r_0 $ the straight up top is unstable.  As $r$ increases past
  $r_0 = \sqrt{-W'(0)}$ the straight up top gains stability.
  Additionally a branch of unstable relative equilibria bifurcates off
  the straight up position.
  \begin{tikzpicture}
    \draw[thick, dashed] (0,0) plot[domain=0:1.4] (\x,{1+(\x)^2}) node[anchor=south west] {};
 \draw[->] (0,0) -- (3,0) node[right] {$u$};
 \draw[thick, dashed] (0,0) -- (0,1) node[left]{$r_0$};   
 \draw[thick,->] (0,1) -- (0,3) node[above] {$r$};
\end{tikzpicture}
\end{itemize}
\end{theorem}
\mbox{}\\
\begin{example}[Lagrange Top]  For the Lagrange top the potential $V\in
  \cin(SO(3))^{SO(2)\times SO(2)}$ is given by $V(A) = \langle Ae_3,
  e_3 \rangle $ for all $A\in SO(3)$.  Hence the corresponding
function $V\in \cin(S^2)^{SO(2)}$ is
\[
V(q) = \langle q, e_3 \rangle = q_3.
\]
Since $V(u, 0, \sqrt{1-u^2}) = \sqrt{1-u^2}$,  the 
function $W\in \cin([0,1))$ is given by
\[
W(t) = \sqrt{1-t}.
\]
Since %
$W'(0)= -\frac{1}{2}<0$ and since $W'(0) < -\frac{1}{4} = W''(0)$
we are in the case (ii) of Theorem~\ref{thm:1.1}: as the rate of spin decreases
the straight up top loses stabilty at $r_0 = 2$ and a stable
relative periodic orbit appears nearby, i.e., the tip of the top will
trace out a circle around the vertical axis.
\end{example}

\begin{example}[Kirchhoff top]   A Kirchhoff top is a family of
symmetric tops with the potential $V(A) = \langle A e_1, e_1\rangle
  ^2 + \langle A e_2, e_2\rangle^2 + c\langle A e_3, e_3 \rangle ^2$
  where $c>0$ is a constant.  It corresponds to the $SO(2)$-invariant
  function $V:S^2\to \R$ given by
  \[
    V(q) = q_1^2 + q_2^2 + c q_3^2.
  \]
Since $V(u, 0,\sqrt{1-u^2}) = u^2 + c(1-u^2) = c + (1-c)u^2$ the
function $W\in \cin([0,1))$ is given by $W(t) = c + (1-c)t$.  Then $W'(0) =
1-c$ and $W''(0)=0$.  Therefore if $c<1$ the derivative  $W'(0) >0$ and there is no bifurcation: the
upright top is stable for all values of $r$.   If $c>1$
\[
  W'(0) = 1-c <0 = W''(0),
\]
so a ``Hamiltonian Hopf'' bifurcation takes place.    
\end{example}

It has been argued that Lagrange and Kirchhoff tops undergo
Hamiltonian Hopf bifurcations as the rate of spin decreases past a
critical value: the sleeping top loses stability and a stable periodic
orbit appears nearby \cite{CvdM, G, BZ}.  Since in
symmetric tops these periodic orbits are in fact relative equilibria
it may be better to view the bifurcation as  figure 8 bifurcations.
These figure 8 
bifurcations are typical for one degree of freedom Hamiltonian systems
with $\Z/2$ symmetries \cite{GoSt}. 
\\
\subsection*{Organization of the paper}
In Section~\ref{sec:2} we review a result of Satzer and Kummer on the
symplectic quotients of $T^*P$ where $S^1\to P\to B$ is a principal
$S^1$ bundle.  We use it reduce the study of symmetric tops to the
study of families of Hamiltonian systems on a magnetic 2-sphere.  This
is an old idea and I include it to keep the paper self-contained.
Section~\ref{sec:3} recalls the theory of singular symplectic
reduction.  We start with the developments in late 1980's -- early 1990's
and continue with more recent work. We then recall the notion of a $\cin$-ring and formulated
the differential spaces of Sikorski in terms of $\cin$-rings.  We 
describe \'Sniatycki's view of singular symplectic quotients as
differential spaces.  In Section~\ref{sec:4} we reduce the study of
families of $S^1$-symmetric Hamiltonian systems on a magnetic 2-sphere
to one variable
calculus. Section~\ref{sec:5} proves the main result of the paper.\\
\subsection*{Acknowledgements} 
I  thank Victor Guillemin for getting me interested in tops and
for making available to me his unpublished notes on the Lagrange top  \cite{G}.  I 
 thank Iliya Zakharevich for many helpful discussions  and 
Richard Cushman for the encouragement.

\section{Regular $S^1$ reduction} \label{sec:2}

It is  well known  that symplectic quotients of the cotangent
bundle  $T^*P$ of a principal bundle $G\to P\to B$ by the lifted action of $G$ 
are symplectic fiber bundles over $T^*B$ with coadjoint orbits of $G$ as
fibers.   If the group in question is 
abelian the proof is simpler since the fibers are points.  The result is due to several people.
 I believe  that the version below,  Theorem~\ref{thm:2.1}, is mostly due to Satzer \cite{Sa} and
Kummer \cite{K}.  Kummer, in turn, relies on the work of Sternberg
\cite{Sternberg} on minimal coupling as reformulated by Weinstein in
\cite{Weinstein}.  I learned the forumalation and the proof of the theorem
 from Victor Guillemin in the mid 1980s.

\begin{theorem} \label{thm:2.1}
Let $P$ be a manifold with a free $S^1$ action, $h(q,p) = \frac{1}{2}
g^*_q (p,p) +V(q)$ an $S^1$-invariant Hamiltonian on the cotangent
bundle $T^*P$ (so $g$ is an $S^1$-invariant metric on $P$ and $V\in
\cin(P)^{S^1}$ an invariant potential) and $\mu:T^*P\to \R =
Lie(S^1)^*$ the associated invariant moment map.   View $P$ as a
principle $S^1$-bundle over $B:= P/S^1$ with the projection $\pi:P\to
B$.

Then for any $r\in \R$ metric $g$ induces a
diffeomorphism
\[
\varphi_r: T^*B \to T^*P/\!/_r S^1 := \mu\inv (r)/S^1
\]
between the cotangent bundle of the base $B$ and the symplectic
quotient.  Moreover
\[
{\varphi_r}^*\omega_r = \omega_{T^*B} + r \tau^*F
\]
where $\omega_r$ is the reduced symplectic form on the quotient
$ T^*P/\!/_r S^1$, $F\in \Omega^2(B)$ is the curvature of the connection 1-form $A$
on $P$ defined by the metric $g$ and $\tau:T^*B\to B$ the canonical
projection.  Finally the pullback by $\varphi_r$ of the induced
Hamiltonian $h_r$ is
\begin{equation} \label{eq:2.2}
({\varphi_r}^*h_r) \,(b, \eta)  = \frac{1}{2} \bar{g}^*_b (\eta,
\eta) +\frac{1}{2}  g^* _q (A_q, A_q) + V(b) 
\end{equation}
for all $b\in B, \eta\in T^*_bB$, where $q\in \pi\inv (b)$ is any point
in the fiber of $\pi$ above $b$ and $\bar{g}$ is the metric induced on
$B$ by $g$.  Note that we are identifying $V\in \cin(P)^{S^1}$ with
the corresponding function on $B$.
\end{theorem}  

\begin{proof}[Sketch of proof]
The zero level set of $\mu:T^*P \to \R$ is the annihilator
$\scV^\circ$ of the vertical bundle of $\pi:P\to B$.   It is not hard
to show that the map
\[
\psi_r: \scV^\circ \to \mu\inv (r),\qquad \psi_r(q, p) = (q, p+r A_q)
 \] 
 is an $S^1$-equivariant diffeomorhism. Here  as above $q\in P$ is a point,
 $p\in \scV^\circ _q$ a covector.   Hence $\psi_r$ descends to a
 diffeomorphism
 \[
\varphi_r: \scV^\circ /S^1 \simeq T^*B \to \mu\inv (r)/S^1. \quad %
\]
The pullback by $\psi_r$ of the tautological 1-form $\theta_{T^*P}$ on
the cotangent bundle of $P$ turn out to be the sum of
$pr^*\theta_{T^*B}$ and $r \,A$ (where $pr:\scV^\circ \to T^*B$ is
the canonical submersion  and we identified $A\in \Omega^1(P)$ with
its pullback to $\scV^\circ$):
\[
\psi_r^* \theta_{T^*P} = pr^* \theta_{T^*B} + r \, A,
\]
see \cite{K}.   Consequently
\[
\varphi_r^* \omega_r  = d\theta_{T^*B} + r dA = \omega_{T^*B} +r \tau^*F.
\]
Finally equation \eqref{eq:2.2} follows from the definition of the
diffeomorhism $\psi_r$ and the fact that the connection 1-form $A$ is
induced by the metric $g$.
\end{proof}

\begin{corollary} \label{cor:2.2}
The symplectic quotient at
$r\in \R = (Lie(SO(2))^*$ of a symmetric top
\[
  (T^*SO(3),
  \omega_{T^*SO(3)}, h(q,p) = \frac{1}{2} g^*_q(p,p) +V(q))\]
with
respect to the action of $SO(2)$ by multiplication on the right is the
classical mechanical system 
\begin{equation} \label{eq:2.4}
(T^*S^2, \omega_{T^*S^2} +r \omega_{S^2}, h(q,p) = \frac{1}{2}
\bar{g}^*_q(p,p) +V(q))
\end{equation}
with $SO(2)$ symmetry, where $\bar{g}$ is an $SO(3)$ invariant metric
on the sphere $S^2$ induced by $g$, $\bar{g}^*$ the dual metric
and $V\in \cin(SO(3))^{SO(2)\times SO(2)}$ is identified with an
$SO(2)$ invariant function on $S^2$ which we again call $V$.
 \end{corollary} 

\begin{proof}
We apply Theorem~\ref{thm:2.1} to the action of $S^1 = SO(2)$ on $P= SO(3)$ by
multiplication on the right.  Then $B = P/S^1$ is the standard
2-sphere $S^2$, the metric $g$ is left $SO(3)$ and right
$SO(2)$-invariant and $V\in \cin(SO(3)) ^{SO(2)\times SO(2)} =
\cin(SO(3)/SO(2))^{SO(2)}$.  Then the induced metric $\bar{g}$ on $S^2$
is $SO(3)$-invariant hence (up to a scalar multiple that we
ignore) is the standard round metric on $S^2$.   The connection 1-form
$A$ and its curvature $F\in \Omega^2(S^2)$ are both
$SO(3)$-invariant.  Hence $F$ is the standard area form on $S^2$
(possibly up to a factor that depends on normalization that we again
ignore).    Finally the function 
\[
q\mapsto \frac{1}{2}g^* (A_q, A_q)
\]
on $SO(3)$ is $SO(3)$-invariant, hence constant.  There is no harm in
dropping it.   We conclude that
under the diffeomorhism $\varphi_r:T^*S^2\to T^*SO(3)/\!/_r S^1$ the
reduced classical mechanical system is 
\[
(T^*S^2, \omega_{T^*S^2} +r \omega_{S^2}, h(q,p) = \frac{1}{2}
\bar{g}^*_q(p,p) +V(q)).
\] 
\end{proof}

\begin{remark}
It is convenient to identify the cotangent bundle $T^*S^2$ of the 2-sphere
with a submanifold of $\R^6$:
\[
T^*S^2 =  \{(q,p) \in \R^3 \times \R^3 \mid \langle q, q\rangle =1,
\langle q, p\rangle =0\} .
\]  
where $\langle \cdot, \cdot\rangle$ is the standard inner product on
$\R^3$.   With this identification the action of $SO(2) $ on $T^*S^2$
becomes identified with the restriction  the diagonal action on $\R^3
\times \R^3$ by rotations about $e_3 = (0,0,1)$.   The corresponding
moment map $\mu_r: (T^*S^2, \omega_{T^*S^2} +r \omega_{S^2}) \to \R$
is given by
\begin{equation} \label{eq:2.3}
\mu_r (q,p) = \langle q\times p, e_3\rangle + r q_3  
\end{equation}
where $\times $ here is the cross product.
\end{remark}

\section{Digression: $\cin$-rings and differential
  spaces} \label{sec:3}

In the next section we will use symplectic reduction of the $SO(2)$-symmetric
Hamiltonian system
\[
  (T^*S^2, \omega_{T^*S^2} +r \omega_{S^2}, \mu_r: T^*S^2 \to \R,
  h(q,p) = \frac{1}{2} g^*_q(p,p) +V(q))
\]
at $r = \mu_r (e_3, 0)$ to
analyze stability and bifurcation of the fixed point $(e_3, 0)$.
Since $r$ is a singular value of the $SO(2)$ moment map $\mu_r$ ($\mu_r$
is given by \eqref{eq:2.3}) the
reduced space $T^*S^2/\!/_r SO(2) = \mu_r \inv (r)/SO(2)$ is singular.
We will argue that a neighborhood of the image of $(e_3, 0)$ in
$T^*S^2/\!/_r SO(2) $ is isomorphic to $T^* ((-1,1))/\Z_2$.   To carry
this out we need to explain what we mean by ``isomorphic.''  In order
to do this we will need to recall the notion of a $\cin$-ring and of a
differential space.

To set the stage we briefly recall the theory of singular reduction as it was
developed in late 1980's -- early 1990s.  It has been known since 1970s that if $\alpha$ is
a regular value of a $G$-equivariant moment map $\mu:M\to \fg^*$ that
arises from a proper action of a Lie group $G$ on a symplectic
manifold manifold $(M,\omega)$ and if the action of the stabilizer
$G_\alpha$ of $\alpha$ on the level set $\mu\inv(\alpha)$ is free then
the $\alpha$-level set is a manifold and the restriction of the
symplectic form $\omega|_{\mu\inv (\alpha)}$ descends to a
symplectic form $\omega_\alpha$ on the symplectic quotient
$M/\!/_\alpha G:= \mu\inv (\alpha)/G_\alpha$.  If additionally
$h\in \cin(M)^G$ is an invariant Hamiltonian then its (local) flow on
$M$ is $G$-equivariant and preserves $\mu\inv (\alpha)$ hence induces
a flow on the quotient $M/\!/_\alpha G$.  On the other hand
$h|_{\mu\inv (\alpha)}$ descends to a smooth function
$h_\alpha\in \cin ( M/\!/_\alpha G)$ and the flow of $h_\alpha$ on
$(M/\!/_\alpha G, \omega_\alpha)$ agrees with the flow induced by $h$.
If the action of $G_\alpha$ on $\mu\inv (\alpha)$ is only locally free
the the level set $\mu\inv (\alpha)$ is still a manifold, 
the quotient $M/\!/_\alpha G$ is naturally a symplectic orbifold and again the
flow on $M/\!/_\alpha G$ induced by an invariant Hamiltonian $h$ is
the flow of the function $h_\alpha$.

If the action of $G_\alpha$ on $\mu\inv (\alpha)$ is not locally free
then the flow of an invariant Hamiltonian $h\in \cin(M)^G$ still
preserves the level set $\mu\inv (\alpha)$ and induces a flow on the
space $M/\!/_\alpha G= \mu\inv (\alpha)/G_\alpha$.  The set
$M/\!/_\alpha G$ is naturally a topological space (take the subspace
topology on $\mu\inv (\alpha)$ and quotient topology on
$\mu\inv (\alpha)/G_\alpha$).  One refers to $M/\!/_\alpha G$ as the
{\sf reduced space at } $\alpha\in \fg^*$ and as a {\sf (singular)
  syplectic quotient}.\\

\noindent
The singular symplectic
quotients are highly structured:

\begin{enumerate}[leftmargin=0.2in]

\item%
  The  space $ M/\!/_\alpha G$ is a symplectic stratified
  space \cite{SL, BL, LW}.  This means that the topological space
  $ M/\!/_\alpha G$ naturally decomposes into a collection of
  symplectic manifolds (these manifolds are called {\sf symplectic
    strata}) and that the singularities of $ M/\!/_\alpha G$ are tame ---
  see \cite{SL, LW}.  More precisely for any Lie subgroup $H$ of $G$
  the intersection $M_{(H)} \cap \mu\inv (\alpha)$ is a manifold
  ($M_{(H)} $ denotes the subset of points of $M$ whose $G$-stabilizer
  is conjugate to $H$).  The  quotient
\[
  (M/\!/_\alpha G )_{(H)}:= (M_{(H)} \cap \mu\inv (\alpha))/G_\alpha
\]
is a manifold as well and restriction of the symplectic form $\omega$
on $M$ to $M_{(H)} \cap \mu\inv (\alpha)$ descends to a symplectic
form on $ (M/\!/_\alpha G )_{(H)}$.  The manifolds $ (M/\!/_\alpha G
)_{(H)}$ are the symplectic strata of $M/\!/_\alpha G $.
\item%
 \label{rmrk:str2} The quotient map $\cp: \mu\inv (\alpha)\to M/\!/_\alpha G $
  induces an isomorphism  $\cp^*: C^0 (M/\!/_\alpha G ) \to C^0 (\mu\inv
  (\alpha))^{G_\alpha}$.  Consequently the preimage of the restriction
  $\cin(M)^G|_{\mu\inv (\alpha)}$ under $\cp^*$ is a subalgebra of  $C^0
  (M/\!/_\alpha G )$ which is denoted by $\cin(M/\!/_\alpha G )$.
  Since
 \[ 
\cin(M/\!/_\alpha G ) \simeq \cin(M)^G/\scI_\alpha
 \]
where $\scI_\alpha = \{f \in \cin(M)^G \mid f|_{\mu\inv (\alpha)} =
0\}$ and since $\scI_\alpha$ is a {\em Poisson} ideal \cite{ACG} the
algebra $\cin(M/\!/_\alpha G ) $ is naturally a Poisson algebra.

\item%
 \label{rmrk:str3} 
  The Poisson bracket on $\cin(M/\!/_\alpha G ) $ is
  compatible with the symplectic forms on the strata of the reduced
  space $M/\!/_\alpha G$: the restriction map
\[
\cin(M/\!/_\alpha G)\to \cin( (M/\!/_\alpha G )_{(H)})
\] 
is Poisson for every subgroup $H<G$ (see \cite{SL}).

\item%
  The compatibility of the Poisson brackets and of the
  symplectic forms on the starta of $M/\!/_\alpha G $ gives two complementary
  ways to view the flow  on the symplectic quotient
  $M/\!/_\alpha G $ induced by an invariant Hamiltonian $h\in \cin(M)^G$.
  One can view it as a collection of Hamiltonian flows on the strata.
  Alternatively the Hamiltonian $h_\alpha\in \cin(M/\!/G)$ induced by
  $h$ defines a derivation
  \[
    \{h_\alpha, -\}: \cin(M/\!/G) \to \cin(M/\!/G),
   \] 
  and the flowlines $\gamma (t)$ of the induced flow are integral
  curves of this derivation:
  \[
\left. \frac{d}{dt}\right|_t f\circ \gamma = \{ h_\alpha, f\} (\gamma (t))
\]
for all functions $f\in \cin(M/\!/G)$.
\end{enumerate}

While the theory described above caused a bit of excitement in 1990s,
a particularly nagging question remained: what were singular
symplectic spaces an example of? %
Several decades later it seems to me that the best answer so far was
found by Jedrzej \'{S}niyaticki \cite{Sn}: singular symplectic
quotients are differential spaces in the sense of Sikorski or, more
generally, $\cin$-locally ringed spaces \cite{Joy}.  It is not a
complete answer and more work remains to be done: see
Remark~\ref{rmrk:more} below.  One can also plausibly argue that
derived differential-geometric symplectic stacks would be another
promising home for singular simplectic quotients.  However at the
moment derived differential geometry is not developed enough for the
study of dynamics on singular symplectic quotients.\\

To explain what differential spaces are 
we start by recalling the notion of a
$\cin$-ring.  The definition below is not standard, but it's
equivalent to the standard one \cite{Joy} and is easier to make sense
of in the first pass (unless you have some background in categorical
universal algebra which I don't).

\begin{definition}\label{def:cring1}
A {\sf $C^\infty$-ring} is a {\em set} $\scC$ together with an
infinite collection of operations 
\[
\left\{g_\scC:\scC^m\to \scC\right\}_{m\geq 0, g\in \cin(\R^m)}
\]
(where
$\scC^0:= \{*\}$, a one-point set and $\R^0: = \{0\}$) so that for all
$n,m\in \N$, all $g\in C^\infty(\R^m)$ and all
$f_1, \ldots, f_m\in C^\infty (\R^n)$
\[
({g\circ(f_1,\ldots, f_m)})_\scC (c_1,\ldots, c_n) =
g_\scC({(f_1)}_\scC(c_1,\ldots, c_n), \ldots, {(f_m)}_\scC(c_1,\ldots, c_n))
\]
for all $(c_1, \ldots, c_n) \in \scC^n$.
Additionally we require that for every coordinate function $x_j:\R^m\to \R$,
\[
(x_j)_\scC (c_1,\ldots, c_m) = c_j.
\]
\end{definition}

\begin{example}
The algebra of functions $\cin(M)$ on a smooth manifold $M$ is a
$\cin$-ring: for any $n>0$, any function $f\in \cin (\R^n)$ one
defines 
\[
f_{\cin(M)} : (\cin(M))^n \to \cin(M)
\]
by 
\[
f_{\cin(M)} (a_1,\ldots, a_n) := f\circ (a_1,\ldots, a_n)
\] 
for any $n$-tuple of 
functions $a_1,\ldots, a_n \in \cin(M)$.
\end{example}

\begin{example}
The real line $\R$ is a $\cin$ ring since it's the algebra of smooth
functions on a one point manifold $*$.  Explicity, given $f\in
\cin(\R^n)$ the corresponding operation $f_\R: \R^n \to \R$ ``is" the
function $f$:
\[
f_\R (a_1,\ldots, a_n) := f(a_1,\ldots, a_n) \quad (\textrm{ the
  evaluation of $f$ on the $n$-tuple $(a_1,\ldots, a_n)$}.
\]
\end{example}

\begin{remark}
  Any $\cin$-ring $\scA$ has an underlying $\R$-algebra.  This is
  because addition and multiplication functions $f(x,y) = x+y$, $g(x,y) = xy$
  are smooth functions as are multiplications by scalars
  $m_\lambda (x) = \lambda x$ (for all $\lambda \in \R$) hence define
  appropriate binary and unary operations on the set $\scA$ making it
  into an $\R$ algebra.

We will not notationally distinguish between $\cin$-rings and their
underlying $\R$-algebras.
\end{remark}

\begin{definition}\label{def:cin-morph}
  A {\sf morphism} of $\scA \to \scB$ of $\cin$-rings is a map of sets
  $\varphi:\scA\to \scB$ which preserves all the operations: for any $n>0$,
  any $a_1,\ldots, a_n\in \scA$ and any $f\in \cin(\R^n)$
  \[
\varphi (f_\scA(a_1,\ldots, a_n)) = f_\scB (\varphi(a_1),\ldots, \varphi(a_n)).
    \]
  \end{definition}
\begin{definition}
A {\sf point} of a $\cin$-ring $\scA$ is a map of $\cin$-rings $p:\scA
\to \R$.
\end{definition}

\begin{example}
Let $M$ be a manifold, $p\in M$ a point and $ev_p:\cin(M)\to \R$ the
evaluation at $p$: $ev_p(f) = f(p)$.  Then $ev_p$ is a point of the
$\cin$-ring $\cin(M)$.
\end{example}

\begin{definition}  \label{def:ideal}\label{def:loc}An {\sf ideal} in a $\cin$ ring $\scA$ is an ideal
  in the underlying $\R$-algebra.
\end{definition}

The following theorem is useful in defining $\cin$-ring structure on singular
symplectic quotients.

\begin{theorem}\label{thm:3.9}
Let $\scA$ be a $\cin$-ring and $I\subset \scA$ an ideal in the
underlying $\R$-algebra.  Then the
quotient $\R$-algebra $\scA/I$ is naturally a $\cin$-ring: for any $n>0$ and any function $f\in C^\infty (\R^n)$ the map
\[
f_{\scC/I}: (\scA/I)^n \to \scA/I, \qquad f_{\scC/I}(c_1 +I, \ldots,
c_n +I):=  f _\scC (c_1,\ldots, c_n) + I
\]
is well-defined for all $(c_1+I, \ldots, c_n+I)\in (\scA/I)^n$.
\end{theorem}

\begin{proof}
The result is well-known.  See \cite{MR} or \cite{Joy} for a proof.
\end{proof}

\begin{definition} \label{def:pd}
A $\cin$-ring $\scA$ is {\sf point determined} if points separate
elements of $\scA$.  That is, if $a\in\scA$ and $a\not = 0$ then
there is a point $p:\scA \to \R$ so that $p(a)\not = 0$.
\end{definition}

\begin{remark} There are many $\cin$-rings that are not point
  determined.  The simplest example is the quotient ring
  $\cin(\R)/\langle x^2\rangle$ where $\langle x^2\rangle$ is the
  ideal generated by the function $x^2$.  This ring has only one point
  $p$ which is given by
  \[
p( f+\langle x^2 \rangle ) = f(0).
 \]   
 See \cite{MR}.
\end{remark}

\begin{definition}  \label{def:sikorski} A {\sf differential space} (in the sense of
  Sikorski) is a pair $(N, \cin(N))$ where $N$ is a topological space and $\cin(N)$ is a
  set of real-valued function on $N$ subject to the following three
  conditions (the notation $\cin(N)$ is meant to be suggestive of a
  smooth structure on the space $N$):
\begin{enumerate}[leftmargin=0.3in]
  \item \label{def:sikorksi:it1} The topology on $N$ is the smallest topology making every
    function in $\cin(N)$ continuous.
\item \label{def:sikorksi:it2} For any $n>0$, any smooth function $f\in \cin(\R^n)$ and any
  $n$-tuple $a_1,\ldots, a_n\in \cin(N)$, the composite $f\circ
  (a_1,\ldots, a_n)$ is in $\cin(N)$.
\item \label{def:sikorksi:it3} For any open cover $\{U_i\}_{i\in I}$ and any function $g:N\to
  \R$ so that for each $i\in I$ there is $a_i\in \cin(N)$ with $g|_{U_i}
  = a_i |_{U_i}$  the function $g$ is in $\cin(N)$.
\end{enumerate}
\end{definition}

\begin{remark}
 Condition (\ref{def:sikorksi:it1}) amounts to requiring that the sets
 $\{ \{f \not =0 \} \mid f\in \cin(N)\}$ generate the topology on
 $N$. This may be viewed as a $\cin$-analogue of the Zariski topology on
 affine varieties.
 One can also show that (\ref{def:sikorksi:it1}) is equivalent to existence
 of bump functions: for any open set $U\subset N$ and for any point
 $x\in U$ there is a function $f\in \cin(N, [0,1])$ with %
 $\supp f\subset U$ and $f$ identically 1 near $x$.
  
 Condition (\ref{def:sikorksi:it2}) %
 amounts to saying that the $\R$-algebra $\cin(N)$
 is in fact a $\cin$-ring.  For some reason most of papers that deal
 with differential spaces never explicitly mention $\cin$-rings.  Note
 that the $\cin$-ring $\cin(N)$ is point determined since it consists
 of actual functions and for any point $p\in N$ the evaluation map
 $ev_p:\cin(N)\to \R$ is a point of the $\cin$-ring $\cin(N)$.

  The third condition can be interpreted as follows: by restricting
  the functions in $\cin(N)$ to open subsets of $N$ one obtains a
  presheaf.  Denote its sheafication by $\cin_N$.  Condition
  (\ref{def:sikorksi:it3}) then amounts to requiring that the
  $\cin$-ring of global section $\cin_N(N)$ of this sheaf {\sf is} the $\cin$-ring
  $\cin(N)$.   Note that in particular a
  differential spaces is implicitly a $\cin$-ringed space.
\end{remark}

\begin{definition} \label{def:diffsp_mor} A {\sf map} or a {\sf morphism} from a
  differential space  $(M,\cin(M))$  to a differential space $(N,
  \cin(N))$ is a map of underlying sets $\varphi:M\to N$ so that for
  any $f\in \cin(N)$ the composite $f\circ \varphi$ is in $\cin(M)$.
\end{definition}

\begin{notation}
Differential spaces and their morphisms form a category which we
denote by $\DiffSp$.
\end{notation}

The reader unfamiliar with $\cin$-schemes should feel free to ignore
the following remark.  It will play no role in the paper.
\begin{remark}  
It is not too hard to show that the category of differential spaces
embeds into the larger category of $\cin$-ringed spaces.
$\cin$-schemes \cite{Dubuc, Joy}  also embed into the category of $\cin$-ringed spaces.
It is not at all clear which differential spaces are $\cin$-schemes
and conversely.

A possible  exception is formed by  affine schemes that come from finitely
generated and point determined $\cin$-rings.  These are exactly the
differential spaces that are isomorphic to closed subsets of Euclidean
spaces \cite{KL}.  This class includes all second countable Hausdorff manifolds.

Part of the problem of comparing differential
spaces and affine $\cin$-schemes  is that for given a differential
space $(N,\cin(N))$ it is not clear that a point $p:\cin(N)\to
\R$ of the $\cin$-ring $\cin(N)$ has to come from evaluation at a
point $x\in N$.  If $N$ is a second countable manifold then any point  $p:\cin(N)\to
\R$ of the $\cin$-ring $\cin(N)$ does come from an evaluation at some
point $x\in N$ by the famous ``Milnor's exercise.''  
\end{remark}

A vector field $v$ on a manifold $M$ can be defined as a  derivation $v:\cin(M)\to
\cin(M)$ of the $\R$-algebra of smooth functions on $M$ with values in
$\cin(M)$: $v$ is
$\R$-linear and for any two functions $f,g\in \cin(M)$
\[
v(fg) = v(f) g + f v(g).
\]
For $\cin$-ring $\scA$ there is another notion of a derivation of
$\scA$ with values in $\scA$:

\begin{definition} \label{def:cin-der-dumb}
  A {\sf $\cin$-derivation} of a $\cin$ -ring $\scA$ is a map $X:\scA\to \scA$ so that 
for any $n>0$, any $f\in \cin(\R^n)$ and any $a_1,\ldots, a_n\in
\scA$
\[
X(f_\scA(a_1,\ldots,a_n)) = \sum_{i=1}^n (\partial_i
f)_\scA(a_1,\ldots, a_n)\,X(a_i).
\]
\end{definition}
\begin{remark}
One can show that for a large class of $\cin$-rings that includes
point determined $\cin$-rings the two notions of derivations coincide.
See \cite{Yamashita}.  Thus if $(N, \cin(N))$ is a differential space
then an $\R$-algebra derivation $v:\cin(N)\to \cin(N)$ is a
$\cin$-ring derivation.
\end{remark}

We are now in position to define integral curves of derivations on
differential spaces.

\begin{definition}
  Let $v:\cin(M)\to \cin(M)$ be a derivation on a differential space
  $M$.  An {\sf integral curve $\gamma$ of $v$ through a point $p\in
    M$ } is
  either a map $\gamma: \{0\} \to M$ with $\gamma(0) = p$ or a smooth map 
  $\gamma: (J, \cin(J))\to (M, \cin(M))$ from an interval
  $J\subset \R$ containing 0  so that
  \begin{equation}
\frac{d}{dx} (f\circ \gamma) = v(f) \circ \gamma
\end{equation}
for all function $f\in \cin(M)$.  Note that unlike \cite{Sn} we do
{\em not} require $J$ to be an open interval.
\end{definition}

\noindent  The reader may be puzzled why we allow integral curves to
only exist for time $t=0$ or to have non-open intervals as domains of
definition.  Example~\ref{ex:4.18} below illustrates why this may be
useful.  Note that for singular symplectic quotients $M/\!/_\alpha G$
and for Hamiltonian $h_\alpha \in \cin(M/\!/_\alpha G)$
induced by $h\in \cin(M)^G$ the maximal  integral curves of the derivation
\[
\{ h_\alpha, \cdot \}: \cin(M/\!/_\alpha G) \to \cin(M/\!/_\alpha G)
\]
are maximal  integral curves of vector fields on manifolds and so have open
intervals as domains of definition.
\begin{example} \label{ex:4.18}
  Let $M$ be the standard closed disk in $\R^2$:
  $M = \{(x,y)\mid x^2 +y^2 \leq 1\}$. Then $M$ is a manifold with
  boundary and  a differential space.  Consider the vector
  field $v = \frac{\partial}{\partial x}$ on $M$.  The integral curve
  of $v$ through $(0,1)$ is $\gamma:\{0\} \to M$, $\gamma(0) = (0,1)$;
  it only exists for zero time.  Note that $v$ does have a flow. It's
  a smooth map $\Phi$ from
  $U =\{((x,y), t) \in \R^2\times \R \mid x^2 + y^2 \leq 1, (x-t)^2
  +y^2 \leq 1\} \to M$.  It is given by $\Phi((x,y),t) = (x+t, y)$.
  Note that here we view $U\subset \R^3$ as a differential space.
  Note also that $U$ is {\sf not} a manifold with corners as one can
  see by looking at its singularities.
\end{example}

As in the case of manifolds related derivations have related integral
curves:
\begin{lemma} \label{lem:3.22}
Let $\varphi:(M,\cin(M))\to (N, \cin(N))$ be a map of differential
spaces, $X:\cin(M)\to \cin(M)$, $Y:\cin(N)\to \cin(N)$ two derivations
so that the diagram
\[
\xy
(-15,10)*+{C^\infty(M)}="1";
(15,10)*+{C^\infty(M)}="2";
(-15,-5)*+{C^\infty(N)}="3";
(15,-5)*+{C^\infty(N)}="4";
{\ar@{->}^{X} "1";"2"};
{\ar@{->}^{\varphi^*} "3";"1"};
{\ar@{->}_{\varphi^*} "4";"2"};
{\ar@{->}^{Y} "3";"4"};
\endxy
\]
commutes.  I.e., for all $f\in \cin(N)$
\[
  Y(f)\circ \varphi = X(f\circ \varphi).
  \]
Then for any integral curve $\gamma:I\to M$ of $X$,
$\varphi\circ \gamma:I\to N$ is an integral curve of $Y$.
\end{lemma}
\begin{proof}
  The proof is the same as in the case of manifolds:
  \[
\frac{d}{dt} (f\circ \varphi \circ \gamma) = X(f\circ \varphi) \circ
\gamma = Y(f)\circ \varphi \circ \gamma.
\]    
\end{proof}  

We now come to the main reason for discussing differential spaces,
which is a theorem due to \'{S}niyticki (see \cite{Sn}, for example).

\begin{theorem}[\'{S}niyticki]
Let $(M, \omega)$ be a symplectic manifold with a proper Hamiltonian
action of a Lie group $G$ and corresponding equivariant moment map
$\mu: M\to \fg^*$.  Then for any $\alpha\in \fg^*$ the quotient
\[
M/\!/_\alpha G:= \mu\inv (\alpha)/G_\alpha
\]
is a Hausdorff differential space with the space of smooth functions
\[
\cin(M \!/_\alpha G):= \cin(M)^G|_{\mu\inv (\alpha)}.
  \]
\end{theorem}

\begin{proof}[Sketch of proof] Since the action of $G$ is proper and
  the stabilizer $G_\alpha$ is closed in $G$, the action of $G_\alpha$
  on $M$ is proper.  Consequently $M/G_\alpha$ and its closed subset
  $\mu\inv (\alpha)/G_\alpha$ are Hausdorff.

  It is easy to check that
  the $\R$-algebra $\cin(M)^G$ of $G$-invariant functions is a $\cin$
  ring: if $a_1,\ldots, a_n \in \cin(M)^G$ are invariant functions and
  $f\in \cin(\R^n)$ then the composite $f\circ (a_1,\ldots, a_n)$ is
  also $G$-invariant.
  Since
\[
  \scI_\alpha: = \{ f\in \cin(M)^G \mid f|_{\mu\inv (\alpha)} \}= 0
\]  
  is an  ideal in $\cin(M)^G$ the quotient $\cin(M)^G/\scI_\alpha$ is
  a $\cin$-ring (see Theorem~\ref{thm:3.9}).  Hence
  \[
    \cin(M)^G|_{\mu\inv (\alpha)} = \cin(M)^G/\scI_\alpha
  \]
  is a $\cin$-ring.

 Conditions  (\ref{def:sikorksi:it1}) and  (\ref{def:sikorksi:it3})
 follow from the existence of a $G$-invariant partition of unity on $M$
 subordinate to a $G$-invariant open cover of $M$.
 \end{proof} 

 \begin{definition}
Two reduced spaces $M/\!/_\alpha G$ and $N/\!/_\beta H$ are {\sf
  isomorphic} if there is an isomorphism
\[
\varphi: (M/\!/_\alpha G, \cin(M/\!/_\alpha G) )\to  (N/\!/_\beta H,
\cin(N/\!/_\beta  H))
\]
of differential spaces (cf.\ Definition~\ref{def:diffsp_mor}) so that
\[
\varphi^*: \cin(N/\!/_\beta  H) \to \cin(M/\!/_\alpha G)
\]
is an isomorphism of Poisson algebras (see \ref{rmrk:str2} ).
\end{definition}   

\begin{lemma} \label{lem:iso-red}
Suppose $
\varphi: (M/\!/_\alpha G, \cin(M/\!/_\alpha G) \to  (N/\!/_\beta H,
\cin(N/\!/_\beta  H))$ is an isomorphism of reduced spaces and $h\in
\cin(N/\!/_\beta  H))$ a Hamiltonian.  Then $\varphi$ sends the
integral curves of $\varphi^*h$ (meaning the integral curves of the
derivation $\{\varphi^*h, \cdot \}:\cin(M/\!/_\alpha G) \to \cin(M/\!/_\alpha G)$)
to the integral curves of $h$.
\end{lemma}

\begin{proof}
The proof follows easily from the definition of an isomorphism of
reduced spaces and Lemma~\ref{lem:3.22}.
\end{proof}

\begin{corollary}\label{cor:3.26}
Suppose $
\varphi: (M/\!/_\alpha G, \cin(M/\!/_\alpha G) \to  (N/\!/_\beta H,
\cin(N/\!/_\beta  H))$ is an isomorphism of reduced spaces and $h\in
\cin(N/\!/_\beta  H))$ a Hamiltonian as in Lemma~\ref{lem:iso-red}
above. Then
\begin{itemize}[leftmargin=0.1in]
  \item A point $x\in N/\!/_\beta H,
\cin(N/\!/_\beta  H)$ is a stable equilibrium of $h$ if and only if
$\varphi\inv (x)$ is a stable equilibrium of $\varphi^*h$.
\item A point $x\in N/\!/_\beta H,
\cin(N/\!/_\beta  H)$ is an unstable equilibrium of $h$ if and only if
$\varphi\inv (x)$ is an unstable equilibrium of $\varphi^*h$.
\end{itemize}   
\end{corollary}
The final bit of theory that we'll need to analyze the stability of
tops is a  theorem Montaldi \cite{Mon}:
\begin{theorem}\label{thm:Montaldi}
  Let $(M, \omega)$ be a symplectic manifold with a Hamiltonian action
  of a compact Lie group $G$ and corresponding equivariant moment map
  $\mu: M\to \fg^*$.  Let $h\in \cin(M)^G$ be an invariant
  Hamiltonian, $\alpha\in\fg^*$,
  $pr: \mu\inv (\alpha) \to M/\!/_\alpha G$ the quotient map and
  $h_\alpha\in \cin(M/\!/_\alpha G)$ the reduced Hamiltonian (so
  $pr^*h_\alpha = h|_{\mu\inv (\alpha)}$.)  Suppose
  $x\in \mu\inv (\alpha)$ is a point so that $pr(x)$ is a local
  minimum or a local maximum of $h_\alpha$.  Then $x$ is relative
  equilibrium of $h$ which is a $G$-stable in $M$.
\end{theorem}

We end the section with a parenthetic remark on differential spaces
and singular reduction.  The remark will play no role in the rest of
the paper.
\begin{remark} \label{rmrk:more}
The differential space approach to singular reduction is not a
complete answer for the following annoying reason.

By a theorem of
Dubuc and Kock \cite{DK} for any $\cin$-ring $\cA$ there exists a
module of $\cin$-K\"ahler differentials $\Omega_\cA^1$ together with a
universal derivation $d_\cA:\cA \to \Omega_\cA^1$: for any
$\cA$-module $\cM$ and any derivation $X:\cA\to \cM$ there exists a
unique $\cA$-module map $\varphi_X:\Omega^1_\cA\to \cM$ so that $X =
\varphi_X \circ d_\cA$.  The universal derivation $d_\cA$ is functorial
in $\cA$.   One can further mimic Grothendieck algebraic de Rham
complex and produce an $\cin$-algebraic de Rham complex
\[
 \Omega^\bullet_\cA\xrightarrow{d} \Omega^{\bullet+1}_\cA,
 \] 
see \cite{KL}. %
Moreover if $\cA
= \cin(M)$ for a manifold $M$ then $\Omega^\bullet_{\cin(M)}$ is the
usual de Rham complex (this is not obvious). Since the complex $\Omega_\cA^\bullet$ is also
functorial in $\cA$,    for any closed subset  $Z$  of a manifold $M$ the
the surjective restriction map
\[
\cin(M)\to \cin(Z), \quad f\mapsto f|_Z
\]
extends to a surjective map $\Omega(M)^\bullet \to
\Omega_{\cin(Z)}^\bullet$ of differential graded algebras.   In
particular if $\mu :(M, \omega) \to \fg^*$ is an equivariant moment
map for a Hamiltonian action of a Lie group $G$ on a symplectic
manifold then the restriction of the 2-form $\omega$ to the zero level
set $\mu\inv (0)$ makes a perfectly good sense: $\omega|_Z$ is a
closed 2-form in the $\cin$-algebraic de Rham complex of the
$\cin$-ring $\cin(\mu\inv (0))$.   This is true regardless of whether  or not
$\mu\inv (0)$ is a manifold.  Furthermore, as we have seen, the
quotient $\mu\inv (0)/G$ is a $\cin$-ring and the quotient map
$\pi: \mu\inv (0)\to \mu\inv (0)/G$ is a smooth map of differential
spaces.   The trouble comes from the fact that unlike the regular case
the 2-form $\omega|_{\mu\inv (0)}$ need not descend to a 2-form on the
quotient $\mu\inv (0)/ G$: there need not exist any 2-form $\sigma \in
\Omega^2_{\cin (\mu\inv (0)/ G)}$ with
\[
\pi^* \sigma = \omega|_{\mu\inv (0)}.
\]  

Here is a simple example: let $(M, \omega) = (\R^2, dx \wedge dy)$, $G
= \{\pm 1\}$ acting by $(-1)\cdot (x,y) = (-x, -y)$.   Then the moment
map $\mu$ is identically 0 and $\mu\inv (0) /G = \R^2/\{\pm 1\}$.
The $\cin$-ring  of smooth functions on the quotient  ``is'' the ring of invariant functions
$\cin(\R^2)^{\{\pm 1\}}$.  With this identification the pullback
$\pi^*: \cin(\R^2/G) \to \cin (\R^2)$ is simply the inclusion
$\cin(\R^2)^{\{\pm 1\}} \hookrightarrow \cin(\R^2)$.  By a theorem of
G.\ Schwarz the $\cin$-ring  of invariant functions is  generated by $x^2, xy$ and
$y^2$.  
Consequently the module of 1-forms $\Omega^1_{\cin(\R^2)^{\{\pm 1\}}}$
on the quotient is generated
by $xdx, ydy$ and $xdy+ ydx$.   It follows that there is no
2-form  $\sigma\in \Omega^2_{\cin(\R^2)^{\{\pm 1\})}}$ with $\pi^*\sigma = dx
\wedge dy$.

There are several ways to fix the problem.  For example on can replace
the ``coarse'' quotient $\R^2/G$ by the stack quotient $[\R^2/G]$.
Then the $G$-invariant 2-form $dx \wedge dy $ does descend to a closed
2-form on the stack quotient.

The example suggests to me that in order to fully understand singular
symplectic quotients one would need understand Hamiltonian dynamics on
stacks. For Deling-Mumford stacks over a site of manifolds this has
been done.  But one would need to consider Artin stacks over  a site
of differential spaces or maybe over a site of $\cin$-affine schemes.
\end{remark}  

\section{Reduction to one variable calculus} \label{sec:4}

We are now back to studying the family \eqref{eq:2.4} of
$SO(2)$-symmetric Hamiltonian systems on the cotangent bundle $T^*S^2$
of the 2-sphere.
Since we are interested in the behavior of straight up tops, we may
restrict our attention to the upper hemisphere
\[
S^2_+:= \{(q_1, q_2, q_3)\in \R^3 \mid q_3>0\}.
\]
The projection
\begin{equation} \label{eq:coord}
  \psi: S^2_+ \to \R^2, \qquad\psi(q) =  (q_1, q_2)
\end{equation}
defines
a coordinate chart on $S^2$.  Denote the corresponding coordinates on
$T^*S^2_+$ by $(x,y, p_x, p_y)$.   In these coordinates the symplectic
form $\omega_r$ is given by
\[
\omega_r = dx\wedge dp_x + dy\wedge dp_y + \frac{r}{\sqrt{1-x^2-y^2}}
  dx \wedge dy,
\]  
 the $SO(2)$ moment map $\mu_r$ is
\[
\mu_r (x,y,p_x, p_y) = xp_y - yp_x +r \sqrt{1-x^2-y^2},
\]
 the metric $g$ is given by the matrix 
\[
g=\left( \begin{array}{lr}
			1+\frac{x^2}{z^2} & \frac{xy}{z^2}\\
			\frac{xy}{z^2}   & 1+\frac{y^2}{z^2}
							\end{array}\right)
\]
and the dual metric $g^*$ by the  inverse  matrix 
\[
g^{-1}= z^2
			\left( \begin{array}{lr}
			1+\frac{y^2}{z^2} & -\frac{xy}{z^2}\\
			-\frac{xy}{z^2}   & 1+\frac{x^2}{z^2}
\end{array}\right).
\]
Here and below
\[
  z= \sqrt{1-x^2 -y^2}.
\]
Consequently the Hamiltonian $h$ of \eqref{eq:2.4} in these
coordinates is
\[
  h(x,y,p_x,p_y) = \frac{1}{2} \left( (1-x^2)p_x^2 - 2xyp_xp_y +
    (1-y^2) p_y^2
    \right) + V(x,y, \sqrt{1-x^2-y^2}).
\]

Equivalently  the diffeomorphism
\[
  \varphi: D^2= \{(x,y)\in \R^2 \mid x^2 +y^2 <1\} \to
  S^2_+, \qquad \varphi(x,y) = (x,y, \sqrt{1-x^2-y^2})
\]
induces a
diffeomorphism $\Phi :T^*D^2\to T^*S_+^2$ of the cotangent bundles,
\[
\Phi^*\omega_r =    dx\wedge dp_x + dy\wedge dp_y + \frac{r}{\sqrt{1-x^2-y^2}}
  dx \wedge dy
\]
while
\[
(\mu_r \circ \Phi) \,(x,y,p_x, p_y) = xp_y - yp_x +r \sqrt{1-x^2-y^2}
\]
and
\begin{equation} \label{eq:4.1}
  h \circ \Phi\, (x,y,p_x,p_y) = \frac{1}{2} \left( (1-x^2)p_x^2 - 2xyp_xp_y +
    (1-y^2) p_y^2
    \right) + V(x,y, \sqrt{1-x^2-y^2}).
\end{equation}

\begin{lemma} \label{lem:iso-of_red}
Consider the map $f:T^* (-1, 1)\to %
D^2\times \R^2\subset
\R^4$ given by
\[
f(u,p_u) = (u,0,p_u, \frac{r}{u}(1- \sqrt{1-u^2})
\]  
 for $u\not =0$ and by $f(0, p_u) = (0,0,p_u, 0)$.  Then $f$ is $\cin$,
$f(T^*(-1,1)) \subset \mu_r\inv (r)$, and $f$ induces an isomorphism
\[
f^*: \cin(T^*S^2_+)^{SO(2)} |_{\mu_r\inv (r)} \to \cin(T^*(-1,1))^{\Z/2}
\]
of $\cin$-rings and of Poisson algebras.
Here and below $\Z/2 =\{\pm 1\}$ acts on $(-1, 1)$ by multiplication
and by the lifted action on $T^*(-1,1)$.
Finally the smooth function $f:T^*(-1,1) \to \mu_r\inv (r)$ induces an isomorphism
\begin{equation}
\bar{f}: T^*(-1,1)/(\Z/2) \to \mu_r\inv (r)/SO(2) = T^*S_+^2 /\!/_r SO(2)
 \end{equation}
 of symplectic reduced spaces.  (We view $T^*(-1,1)/(\Z/2)$ as the
 reduction at zero of the cotangent bundle
 $(T^*(-1,1), \omega_{T^*(-1,1)})$ by the lifted action of $\Z/2$.)
\end{lemma}

\begin{proof}
Since $\sqrt{1+x}$ is analytic for $|x|<1$ and since $\sqrt{1+x} = 1
+\frac{1}{2} x + h.o.t.$, the function $u\mapsto
\frac{1-\sqrt{1-u^2}}{u}$ is analytic for $|u|<1$.  Hence $f$ is
$\cin$.

Given $(x,y,p_x, p_y) \in T^*S_+^2$ there is $A\in SO(2)$ so that
$A\cdot (x,y,p_x, p_y) = (x', 0, p'_x,p'_y)$ for some $x',p_x',
p_y'$ with $x'$ unique up to sign. Also if $(x,y,p_x, p_y) \in \mu_r\inv (r)$ and $y=0$ then
$xp_y+r\sqrt{1-x^2} = r$ hence $p_y= r\frac{1-\sqrt{1-x^2}}{x}$.   It
follows that $f(T^*(-1,1)) \subset \mu_r\inv (r)$ and that the $SO(2)$
orbits in $\mu_r\inv (r)\smallsetminus \{(0,0,0,0)\}$ intersect the
image of $f$ in exactly two points.  %
Hence the 
composite map $\pi\circ f: T^*(-1,1) \to \mu_r\inv (r)/SO(2)$ (where
$\pi: \mu_r\inv (r)\to \mu_r\inv (r)/SO(2)$ is the orbit map) descends
to a continuous bijection $\bar{f}: T^*(-1,1)/(\Z/2)\to \mu_r\inv
(r)/SO(2)$.

Since for any $h\in \cin(T^*S^2_+)^{SO(2)} $ the pullback
$f^*h \in \cin(T^*(-1,1))^{\Z/2}$ the map $\bar{f}$ is a map of
differential spaces.  Since for any $(u,p_u)\in T^*(-1,1)$
\[
f( (\Z/2)\cdot (u,p_u))  = \left(SO(2)\cdot f(u,p_u) \right)\cap \mu_r\inv (r)
\]
the pullback map $f^*: \cin(T^*S^2_+)^{SO(2)} |_{\mu_r\inv (r)} \to
\cin(T^*(-1,1))^{\Z/2}$ is injective.   It remains to prove that $f^*$
is bijective and preserves the Poisson brackets.

The preservation of brackets follows from the fact that
\[
f^*\left( dx\wedge dp_x + dy\wedge dp_y + \frac{r}{\sqrt{1-x^2-y^2}}
  dx \wedge dy\right) = du\wedge dp_u
\]  
and (\ref{rmrk:str3}). 

We now argue that $f^* : \cin(T^*D^2)^{SO(2)} \to
\cin(T^*(-1,1))^{\Z/2}$ is onto.  Thanks to Theorem~1 of \cite{Sch} we
know that for a representation $G\to GL(V)$ of a compact Lie group on
a finite dimensional real vector space the $\cin$-ring $\cin(V)^G$ of
invariant functions is finitely generated: there is $k>0$ and
$\sigma_1,\ldots, \sigma_k\in \cin(V)^G$ so that for any $a\in
\cin(V)^G$ there is a smooth function $f\in \cin(\R^n)$ with 
\[
a = f_{\cin(V)^G} (\sigma_1,\ldots, \sigma_k) = f\circ  (\sigma_1,\ldots, \sigma_k) .
\]
Moreover the generators may be taken to be the generators of ring of
invariant polynomials $\R[V]^G$ on $V$.  For the lift of the standard
action of $SO(2)$ on $\R^2$ to  $T^*\R^2 = \R^2\times \R^2$ the $\cin$-ring
$\cin(T^*\R^2)^{SO(2)}$ is generated by four polynomials:
\[
x^2+y^2, \quad p_x^2 +p_y^2, \quad xp_x+ yp_y,\quad  xp_y -y p_x.
\]
For the action of $\Z/2$ on $T^*R = \R^2$ the $\cin$-ring 
$\cin(T^*\R)^{\Z/2}$ is generated by the three polynomials:
\[
u^2, \quad p_u^2, \quad up_u.
\]
Observe that $u^2 = f^* (x^2 +y^2)$, $up_u = f^* (xp_x +yp_y)$ while
$f^*(xp_y- yp_x) = r(1-\sqrt{1-u^2})$.  Hence 
\[
p_u^2 = f^* \left( (p_x^2 + p_y^2) - (xp_y -yp_x)^2(x^2 +y^2) \right).
\]
It follows that 
\[
\left.\cin(T^*\R)^{\Z/2}\right|_{T^*(-1,1)} =f^* \left(
  \left. \cin(T^*\R^2)^{SO(2)}\right|_{\mu_r\inv (r)}\right).
\]
This doesn't quite prove surjectivity of $f^*: \cin(T^*D^2)^{SO(2)} \to
\cin(T^*(-1,1))^{\Z/2}$ since $\cin(T^*(-1,1))^{\Z/2}$ is bigger than
$\left.\cin(T^*\R)^{\Z/2}\right|_{T^*(-1,1)} $.     On the other hand,
since $T^*(-1,1)$ is open in $T^*\R$
the Localization Theorem of Mu\~noz D\'iaz and Ortega
\cite{MO} (see also \cite[p.\ 28]{NGSS}) implies that 
given a function $k\in \cin(T^*(-1,1))^{\Z/2}$ there exist $g, h\in \cin(T^*\R)$ so that
\[
  \left( h|_{T^*(-1,1)}\right) \, k = g|_{T^*(-1,1)} 
\]
and $h|_{T^*(-1,1)} $ is invertible in $\cin({T^*(-1,1)} )$.  By
averaging over $\Z/2$ if neccesary we may assume that $g$ and $h$ are
in $\cin(T^*\R)^{\Z_2}$.  This implies that there are $\tilde{g},
\tilde{h} \in \cin(T^*\R^2)^{SO(2)}|_{\mu_r\inv(r)}$ with
\[
h|_{T^*(-1,1)} = f^* \tilde{h}, \quad g|_{T^*(-1,1)} =f^*\tilde{g}.
\] 
Therefore
\[
 k = \frac{g|_{T^*(-1,1)} }{h|_{T^*(-1,1)}} = f^* \left(\frac{\tilde{g}}{
\tilde{h}}\right) 
\]
and we are done.
\end{proof}

\begin{lemma} \label{lem:4.5}
  If $u\in (-1,1)$ is a critical point of the function
  \[
U_r (u) = 
  \frac{r^2}{2}\left( \frac{1 -\sqrt{1-u^2}}{u} \right)^2 +  V(u,0,\sqrt{1-u^2}).
\]  
(where $r\in \R$ is a parameter) then $((u,0), (0, r \frac{1
  -\sqrt{1-u^2}}{u}))\in T^*S^2_+ \subseteq T^*S^2$ is a relative
equilibrium of the $SO(2)$-invariant Hamiltonian  system
\eqref{eq:2.4} (where we used the coordinates \eqref{eq:coord} on the
upper hemisphere and the induced coordinates on its cotangent bundle).

Moreover if $u$ is a local minimum of $U_r(u)$ then the corresponding
relative equilibrium is relatively stable and if $u$ is a local
maximum of $U_r(u)$  then the correponding relative equilibrium is unstable.
\end{lemma}

\begin{proof}
It follows from Corollary~\ref{cor:3.26} that in order to analyze the stabilty of the
relative equilibria of \eqref{eq:2.4} near the straight up position  it is enough to analyze relative equilibria of
\[
  (T^*(-1,1)/(\Z/2), \cah_r:= f^*(h\circ
  \Phi)\in \cin(T^*(-1,1))^{\Z/2}),
\]
where $h\circ \Phi \in \cin(T^*D^2)^{SO(2)}$ is given by
\eqref{eq:4.1}.  That is, it's enough to analyze the critical points
of the function $\cah_r$ on the manifold $T^*(-1,1)$ (and remember not
to double count since we want to analyze the equilibria on the quotient $T^*(-1,1)/(\Z/2)$).
Since 
\begin{equation}
  \begin{split}
\cah_r (u,p_u) = f^*(h\circ \Phi)\, (u,p_u) = &\frac{1}{2}\left( (1-u^2)p_u^2 +
  r^2\left(\frac{1 - \sqrt{1-u^2}}{u}\right)^2\right) + V(u,
0,\sqrt{1-u^2})\\
=& \frac{1}{2}(1-u^2)p_u^2+ \frac{r^2}{2}\left( \frac{1
    -\sqrt{1-u^2}}{u} \right)^2 +  V(u,0,\sqrt{1-u^2})
\end{split}
 \end{equation}
the critical points of $\cah_r$ are of the form $(u,0)$ where $u$ is a
critical point of the effective potential
\[
U_r (u) = 
  \frac{r^2}{2}\left( \frac{1 -\sqrt{1-u^2}}{u} \right)^2 +  V(u,0,\sqrt{1-u^2}).
\]  
By Montaldi's theorem
(Theorem~\ref{thm:Montaldi}) any extremal points of the Hamiltonian
$\cah_r$ corresond to stable relative equilibria of the system
\eqref{eq:2.4} and any unstable equilibria of $\cah_r$  correspond to
unstable relative equilibria of \eqref{eq:2.2}.

If $u$ is a local  minimum of $U_r$ then $(u,0)$ is a minimum of
$\cah_r$, hence corresponds to a stable relative equilibrium of
\eqref{eq:2.4}.
On the other hand if $u$ is a local  maximum of $U_r$, then $(u,0)$ is
a saddle point of $\cah_r$ and therefore corresponds to an unstable
relative equilibrium of \eqref{eq:2.4}.  
\end{proof}

\begin{remark}
Note that, as we
observed in the introduction, there is a unique function $W\in
\cin([0,1))$ so that
\[
W(u^2) = V(u, 0,\sqrt{1-u^2}).
 \] 
Therefore the effective potential $U_r(u)$ that we need to analyze is
of the form
\begin{equation} \label{eq:4.7}
U_r (u) = 
  \frac{r^2}{2}\left( \frac{1 -\sqrt{1-u^2}}{u} \right)^2 + W(u^2)
\end{equation}
for some function $W\in
\cin([0,1))$  that depends on the top we are studying (cf.\ the
introduction).  We analyze \eqref{eq:4.7} in the next section.
\end{remark}

\section{Analysis of critcal points of $
U_r (u) = 
  \frac{r^2}{2}\left( \frac{1 -\sqrt{1-u^2}}{u} \right)^2 + W(u^2)$
  and a proof of \protect Theorem~\ref{thm:1.1}} \label{sec:5}

\begin{lemma} Let $W\in \cin([0,1))$ be a smooth function and $U_r\in \cin ((-1,1)) $ be given by 
\[
U_r (u) =  \frac{r^2}{2}\left( \frac{1 -\sqrt{1-u^2}}{u} \right)^2 + W(u^2).
\]
Without loss of generality assume that $0\leq r$.
\begin{enumerate}
\item \label{5.1.i}If $W'(0) >0$ then $u=0$ is a local minimum of $U_r$ for all $r$
  and there are no other critical points of $U_r$ near $u=0$.
\item \label{5.1.ii}Suppose $W'(0)<0$ and $W''(0) > W'(0)$.  Then for
  $r>r_0:= \sqrt{-8W'(0)} $ the point $u=0$ is a local minimum of
  $U_r(0)$ and there are no other critical points of $U_r(u)$ for
  $|u|\ll 1$.  For $r<r_0$ the point $u=0$ is a local maximum of the
  function $U_r(u)$ and there is a continuous function $u(r)$ defined
  for $r<r_0$ and $|r-r_0|\ll 1$ so that $u= u(r)$ is a local minimum
  of $U_r(u)$:\\
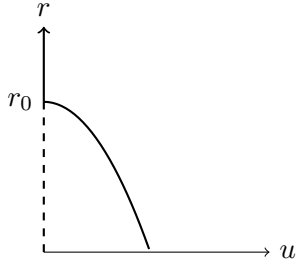
   \begin{tikzpicture}
   \draw[thick] (0,0) plot[domain=0:1.4] (\x,{2-(\x)^2}) node[anchor=south west] {};
\draw[->] (0,0) -- (3,0) node[right] {$u$};
\draw[thick, dashed] (0,0) -- (0,2) node[left]{$r_0$};   
\draw[thick,->] (0,2) -- (0,3) node[above] {$r$};
\end{tikzpicture}
\item \label{5.1.iii}Suppose $W'(0)<0$ and $W''(0) < W'(0)$.  Then  for
  $r<r_0:= \sqrt{-8W'(0)} $ the point $u=0$ is a local maximum of
  $U_r(u)$.  For $r>r_0$ the point $u=0$ is a local minumum and there is a continuous function $u(r)$ defined
  for $r>r_0$ and $|r-r_0|\ll 1$ so that $u= u(r)$ is a local maximum
  of $U_r(u)$:\\
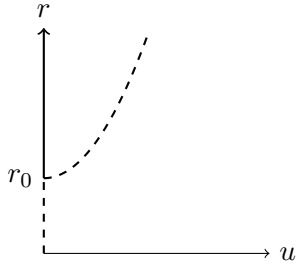
   \begin{tikzpicture}
   \draw[thick, dashed] (0,0) plot[domain=0:1.4] (\x,{1+(\x)^2}) node[anchor=south west] {};
\draw[->] (0,0) -- (3,0) node[right] {$u$};
\draw[thick, dashed] (0,0) -- (0,1) node[left]{$r_0$};   
\draw[thick,->] (0,1) -- (0,3) node[above] {$r$};
\end{tikzpicture}
\end{enumerate}
\end{lemma}

\begin{proof}
The function
\[
v(u)= \begin{cases} \frac{1 -\sqrt{1-u^2}}{u}
    \quad \mathrm{if} \, \, u\not =0 \\
   0 \quad \mathrm{if} \, \, u=0
   \end{cases}
\]
is analytic for $|u|<0$ and is invertible with the inverse $k(v) =
\frac{2v}{1+v^2}$.  Then
\[
U_r(k(v)) = \frac{r^2}{2} v^2 + W(k(v)^2) = \frac{r^2}{2} v^2 + W( \frac{4v^2}{(1+v^2)^2}).
\]  
Let
\[
  f(s) =  W( \frac{4s}{(1+s)^2})
\]
and 
\[
S(r,v) := U_r(k(v)) = \frac{r^2}{2} v^2 + f(v^2).
\]
Then 
\[
f'(0) = 4W'(0) \quad \textrm{and} \quad f''(0) = 16(W''(0) -W'(0)).
\]
Now
\[
  \partial_v S = r^2 v + 2v f'(v^2) = v (r^2 + 2f'(v^2))
\]
and
\[
  \partial_v ^2 S = r^2 + 2f'(v^2) + 4v^2 f''(v^2).
\]
\mbox{}\\
\underline{Suppose that   $f'(0) >0$}. Then $r^2 + 2f'(v^2) >0$ for
all $|v|\ll 1$ and
\[
  \partial_v ^2 S (r,0) = r^2 + 2f'(0) >0.
\]
It follows that  $v=0$ is an isolated local minumum of $S_r(v) = S(r,v)$ for
all $r>0$. Therefore if $W'(0) >0$ then $u=0$ is an isolated  local
minimum of $U_r$ for all $r$.  This proves \eqref{5.1.i}.\\[6pt]
\noindent
\underline{Suppose next that $f'(0) < 0$}.  Then for $r>r_0 = \sqrt{ -2 f'(0)}$ $(=\sqrt{ - 8W'(0)} $)
\[
  \partial_v ^2 S (r,0) = r^2 + 2f'(0) >0
\]
and $r<r_0$
\[
  \partial_v ^2 S (r,0) <0.
\]
It follows that $u=0$ is a local minimum of $U_r(u)$ for
$r> r_0$ %
 and a local maximum for $r<r_0$. %
 We now break up the case $f'(0)<0$ into two generic subcases.
\\[6pt]
\noindent
\underline{Suppose that $f''(0) >0$}.   Then $f'(t)$ is an increasing
function of  $t$ for $|t|\ll 1$.  Consequently $f'(t)$ is invertible
on a neighborhood of $0$
and $(f')\inv $ is also an increasing function in a small neighborhood
of $f'(0) = -r_0^2/2$.  It follows that
\[
\begin{split}
\textrm{ since }-r^2/2 < -r_0^2/2\quad \textrm{ for } r>r_0\geq 0, \qquad &(f')\inv (-r^2/2) < (f')\inv (-r_0^2/2) = 0,
\\
\textrm{ and for } 0 \leq r<r_0\qquad& (f')\inv (-r^2/2) > (f')\inv (-r_0^2/2) = 0.
\end{split}    
\]
Hence for $r>r_0$ the equation
\[
v^2 = (f')\inv (- r^2/2) 
 \]
 has no solutions.  Consequently
 \[
  \partial _v S(v,r) = v(r^2 + 2f'(v^2)) =0 
\]
only if $v=0$.   On the other hand if $r<r_0$ then 
\[
v(r) = \left( (f')\inv (-r^2/2) \right)^{1/2} 
\]
solves
\[
r^2 +2 f'(v^2) =0.
\]
And then
\[
(\partial_v^2 S ) (r, v(r)) = r^2 + 2f' (v(r)^2) + 4 v(r)^2
f''(v(r)^2) =  4 v(r)^2f''(v(r)^2).
\]
Since $f''(v(r_0)^2) = f'' (0) > 0$ by our assumption and since the
function $v(r)$ is
continuous
\[
f''(v(r)^2 ) >0
\]
for $r< r_0$ and $r_0 -r \ll 1$.  Hence $v(r)$ is a local minimum of
$S_r(v) = S(r, v)$.  It follows that $u(r) = \frac{2v(r) }{1+v(r)^2}$
is a local minimum of $U_r$ if $W''(0) > W'(0)$.  This proves \eqref{5.1.ii}. \\[6pt]
\noindent
\underline{Suppose that $f''(0) <0$} (which corresponds to  the case
of $W''(0) < W'(0)$).   Then $f'(t)$ is a decreasing function for $t$
small.  Consequently $f'$ is invertible on a neighborhood of $0$ and
$(f')\inv $ is also an decreasing  function in a small neighborhood
of $f'(0) = -r_0^2/2$.   Then, 
\[
\begin{split}
  \textrm{ since }-r^2/2 < -r_0^2/2\quad \textrm{ for } r>r_0 \geq 0, \qquad
  &(f')\inv (-r^2/2) >(f')\inv (-r_0^2/2) = 0,
\\
\textrm{ and for } r<r_0\qquad& (f')\inv (-r^2/2) < (f')\inv (-r_0^2/2) = 0.
\end{split}    
\]
Hence for $r<r_0$ the equation
\[
v^2 = (f')\inv (- r^2/2) 
 \]
 has no solutions.  Consequently
 \[
  \partial _v S(v,r) = v(r^2 + 2f'(v^2)) =0 
\]
only if $v=0$. 
On the other hand if $r<r_0$ then 
\[
v(r) = \left( (f')\inv (-r^2/2) \right)^{1/2} 
\]
solves
\[
r^2 +2 f'(v^2) =0.
\]
And then
\[
  (\partial_v^2 S ) (r, v(r)) %
  =  4 v(r)^2f''(v(r)^2).
\]
Since $f''(v(r_0)^2) = f'' (0) > 0$ by our assumption,
\[
f''(v(r)^2 ) >0
\]
for $r> r_0$ and $r -r_0\ll 1$.  Hence $v(r)$ is a local maximum of
$S_r(v) = S(r, v)$.  It follows that $u(r) = \frac{2v(r) }{1+v(r)^2}$
is a local maximum of $U_r$.  This proves \eqref{5.1.iii}. 
\end{proof}

\begin{proof}[Proof of Theorem~\ref{thm:1.1}]
By Lemma~\ref{lem:4.5} local minima of the function $U_r$ correspond
to relatively stable relative equilibria of the system \eqref{eq:*}
and local maxima to unstable relative equilibria of \eqref{eq:*}.  By
\eqref{5.1.i} if $W'(0)>0$ the $u=0$  is a local minimum of $U_r$ for
all $r$.  Hence $((u,0), (0,0))\in T^*D^2  \simeq T^*S^2_+$ is a
stable relative equilibrium of \eqref{eq:*} for all values of $r$.
Similarly \eqref{5.1.ii} translates into case (ii) of the theorem and
\eqref{5.1.iii} translates into case (iii).
\end{proof}

\end{document}